\documentclass{conm-p-l}
\usepackage{amssymb,amscd}

\copyrightinfo{1998}{Jonathan Rosenberg}

\newtheorem{thm}{Theorem}[section]    
\newtheorem{lem}[thm]{Lemma}          
\newtheorem{prop}[thm]{Proposition}

\newtheorem{ques}[thm]{Question}
\theoremstyle{definition}

\newtheorem{exm}[thm]{Example}
\newcommand\bC{{\mathbb C}}
\newcommand\bZ{{\mathbb Z}}
\newcommand\bZt{{\mathbb Z/2}}
\newcommand\bR{{\mathbb R}}

\newcommand\cH{{\mathcal H}}
\newcommand\fp{{\mathfrak p}}
\newcommand\fq{{\mathfrak q}}
\newcommand\co{\colon\,}
\newcommand\Cl{\mathrm{Cliff}_{\bC}}
\newcommand\twedge{\textstyle{\bigwedge}}
\newcommand\longlongarrow{\relbar\joinrel\relbar\joinrel\rightarrow}
\newcommand\Hom{\mathrm{Hom}\,}
\newcommand\Ext{\mathrm{Ext}\,}

\newcommand\Gind{G\hbox{-}\mathrm{Ind}\,}
\newcommand\cotan{T^*\!}

\numberwithin{equation}{section}

\begin{document}

\title{The $G$-Signature Theorem Revisited}

\author{Jonathan Rosenberg}
\address{Department of Mathematics, University of Maryland, 
College Park, MD 20742}
\email{jmr@math.umd.edu}
\urladdr{http://www.math.umd.edu/\raisebox{-.6ex}{\symbol{"7E}}jmr}
\thanks{This research was partially supported by NSF Grant \# DMS-96-25336.}

\dedicatory{Dedicated to Mel Rothenberg on his 65th birthday}
\subjclass{Primary 58G12; Secondary 19K33, 19K35, 19K56, 19L47, 
57R91, 57S17, 57R85}
\keywords{$G$-signature, signature operator, equivariant $K$-homology,
equivariant $KK$-theory, characteristic classes, oriented $G$-bordism}
\begin{abstract}
In a strengthening of the 
$G$-Signature Theorem of Atiyah and  Singer, we compute,
at least in principle (modulo certain torsion of exponent
dividing a power of the order of $G$),
the class in equivariant $K$-homology of the signature operator
on a $G$-manifold, localized at a prime idea of $R(G)$, in terms 
of the classes in non-equivariant $K$-homology of the signature 
operators on fixed sets. The main innovations are that the calculation
takes (at least some) torsion into account, and that we are able to extend
the calculation to some non-smooth actions.
\end{abstract}

\maketitle

\section{Introduction}\label{sec:intro}

Let $G$ be a finite group. In studying actions on $G$ on closed 
manifolds $M^n$,\footnote{For the moment we assume the manifold 
and the action to 
be smooth, though we will relax this later on.}\ one of the most important
tools, which comes from analysis of the signature operator, is the 
$G$-Signature Theorem of Atiyah and  Singer (\cite{AS3}, \S6). Suppose $M$
is oriented and the $G$-action preserves the orientation. If we fix a 
$G$-invariant Riemannian metric on $M$, then we can construct the
signature operator $D_M$ (or simply $D$, if $M$ is understood), a
$G$-invariant elliptic
first-order differential operator acting on $\twedge T^*_\bC M$,
given simply by $d+d^*$ (exterior differentiation plus its adjoint with
respect to the metric) together with a $\bZt$-grading of $\twedge T^*_\bC M$
determined by the Hodge $*$-operator (which in turn depends on the orientation
and the metric). Suppose further that $n$, the dimension of $M$, is even.
Then by the formalism of Kasparov theory (\cite{K} --- see also \cite{Ro1}
for a quick summary and \cite{Bl} and \cite{Hig} for good detailed 
expositions),
$D$ determines an equivariant $K$-homology class $[D]\in K^G_0(M)$
which is independent of the choice of metric. If $E$ is any $G$-vector
bundle on $M$, then a choice of a connection on $E$ enables us to define
the signature operator $D_E$ with coefficients in $E$, which again has a 
class $[D_E]\in K^G_0(M)$ independent of the connection. If $c\co M\to 
\mathrm{point}$ is the ``collapse'' map, then the $G$-index
of $D_E$ is simply $\Gind D_E=
c_*([D_E])\in K^G_0(\mathrm{point})=R(G)$, and in turn
this is just the pairing $\langle [D_E],\,[E]\rangle$, where $[E]\in K^0_G(M)$
is the equivariant $K$-theory class of the bundle $E$. The Atiyah-Singer
$G$-Signature Theorem computes this in terms of characteristic class data
of $M$ and $E$.

However, this is only part of the story. Computing $\Gind D_E$ for all
possible $G$-vector bundles amounts to computing the image of $[D]$
under the index map
\begin{equation}\label{eq:Gindex}
\Gind\co K^G_*(M) \to \Hom_{R(G)}(K_G^{-*}(M),\,R(G)).
\end{equation}
However, as pointed out in \cite{Mad} and \cite {MadRos}
(based in part on unpublished work of M.\  B\"okstedt), the map of
equation (\ref{eq:Gindex}) fits inside a short exact universal coefficient
sequence
\begin{equation}\label{eq:UCT}
0 \to \Ext_{R(G)}^1(K_G^{1-*}(M),\,R(G))
\to K^G_*(M) \stackrel{\Gind}{\longlongarrow} 
\Hom_{R(G)}(K_G^{-*}(M),\,R(G)) \to 0,
\end{equation}
and there are canonical isomorphisms
\begin{equation}\label{eq:changeofring}
\Hom_{R(G)}(N,\,R(G))\cong \Hom_\bZ(N,\,\bZ),\qquad
\Ext_{R(G)}^1(N,\,R(G))\cong \Ext_\bZ^1(N,\,\bZ),
\end{equation}
for any $R(G)$-module $N$. Since $R(G)$ is a finitely generated free
$\bZ$-module, and $K_G^{-*}(M)$ is finitely generated over $R(G)$ (since
$M$ has the $G$-homotopy type of a finite $G$-CW-complex), it follows that
the kernel of the index map $\Gind$ of equation (\ref{eq:Gindex}) is precisely
the $\bZ$-torsion in $K^G_*(M)$. Thus the Atiyah-Singer theorem only
computes $[D]\in K^G_0(M)$ modulo torsion.

However, even in the non-equivariant setting, the torsion part of the
$K$-homolo\-gy class of the signature operator is a very interesting
invariant of a manifold \cite{RW2}. The purpose of this paper is therefore
to get more precise information on the $K$-homology class $[D]\in K^G_0(M)$.

\section{Bordism invariance}\label{sec:bord}
While the signature operator on a manifold $M$ with even dimension $n=2k$
is usually described using the de Rham
complex, to describe the signature operator for manifolds of all dimensions
it is convenient to use an equivalent approach using
Clifford algebras (\cite{LM}, Ch.\ II, Example 6.2). By means of the usual
identification of the exterior algebra and Clifford algebra
(as vector spaces, of course, not as algebras), we can view $D_M$ as 
being given by the Dirac operator on $\Cl \cotan M$, the 
complexified Clifford algebra bundle of the cotangent 
bundle (with connection and metric coming from the Riemannian
connection and metric), with grading operator $\tau$ given by the
``complex volume element'' (\cite{LM}, pp.\ 33--34 and 135--137),
a parallel section of $\Cl \cotan M$
which in local coordinates is given by $i^k e_1\cdots e_n$, where
$e_1,\,\ldots,\, e_n$ are a local orthonormal frame for the
cotangent bundle. When the dimension $n=2k+1$ of $M$ is odd, 
$\tau=i^{k+1} e_1\cdots e_n$ 
acting on $\Cl \cotan M$ by Clifford multiplication still satisfies $\tau^2=1$,
but the Dirac operator commutes with $\tau$. Furthermore, if $\sigma$
is the usual grading operator on $\Cl \cotan M$ (which is $(-1)^p$ on products
$e_{i_1}\cdots e_{i_p}$), then $\tau$ and the Dirac operator both
anticommute with $\sigma$. So in this case
we consider the Dirac operator on $\Cl \cotan M$, 
with the grading given by $\sigma$, but with the extra action of the
Clifford algebra $C_1=\Cl \bR$, where the odd generator of $C_1$ acts 
by $\tau$. So for $n$ odd, the signature operator gives a class in the
$K$-homology group $K_1(M)$. (If $M$ is non-compact and we use
a complete Riemannian metric, the $K$-homology class of the signature
is still well-defined and independent of the metric, but lives in 
\emph{locally finite} $K$-homology.)

Now suppose that a finite group $G$ acts on $M^n$, preserving the
orientation. If we choose a $G$-invariant Riemannian metric on $M$,
the signature operator becomes $G$-equivariant, and defines a class
$[D_M]$ living in $K_0^G(M)$ if $n$ is even, $K_1^G(M)$ if $n$ is
odd. In this paper we will be interested in computing $[D_M]$ as
precisely as possible, including torsion information.

An important fact in this context, observed for example in
\cite{PeRW} or in \cite{RW2}, is that if $M^n$ is the boundary
of a compact manifold with boundary $W^{n+1}$, then if $n$ is even,
$[D_M]$ is the image of $[D_{\text{int}\,W}]$ under the boundary map
$K_1^G(W,\partial W)\to K_0^G(\partial W=M)$. However, if $n$ is odd,
the image of $[D_{\text{int}\,W}]$ under the boundary map
$K_0^G(W,\partial W)\to K_1^G(\partial W=M)$ is {\em twice\/} $[D_M]$.
Nevertheless, $[D_M]$ will in this case be the boundary of an operator
in $K_0^G(W,\partial W)$ if $W$ admits a $G$-invariant nonvanishing
vector field pointing, say, inward on $\partial W=M$, or in other
words if $M$ bounds in the sense of equivariant
\emph{Reinhart bordism} as studied in \cite{Kom} and \cite{WW2}. 
(The point is that the
vector field can be used to get a further splitting of $\Cl \cotan W$.)
So we can summarize this information in the following:
\begin{thm}[Bordism Invariance, cf.\ \cite{RW2}]\label{bordinv}
$\!\!$Let $M^n$ be a closed oriented manifold, and suppose a finite group
$G$ acts smoothly on $M$, preserving the
orientation. Suppose given a map $f\co M\to X$. Then if $n$ is even,
$f_*([D_M])\in K_0^G(X)$ only depends on the bordism class of $f$
in $\Omega^G_n(X)$. 
\end{thm}
\begin{proof}It clearly suffices to show that $f_*([D_M])=0$ when
$M=\partial W$ and $f$ extends to a map $g\co W\to X$. 
We use the commutative diagram
\[
\begin{CD}
K_{1}^G{(W,\,\partial W)} @>{g_*}>> K_{1}^G{(X,\,X)}=0\\
@V{\partial}VV @V{\partial}VV\\
K_0^G(M) @>{f_*}>> K_0^G(X)
\end{CD}
\]
and observe that $[D_{\text{int}\,W}]$ maps to $f_*([D_M])$
going down and then across, and to $0$ going across and then down.
\end{proof}
If $n$ is odd, the proof of Theorem \ref{bordinv} only gives
bordism invariance up to a factor of $2$. However, we assume
in addition that $G$ is
of odd order, then in some cases the results of \cite{Kom} and
\cite{WW2} can be used replace oriented $G$-bordism by
oriented Reinhart $G$-bordism, and then the argument will go
through.

We should point out that using Theorem \ref{bordinv}, one can
get an approach to a refinement of the Atiyah-Singer Theorem
using localization in equivariant bordism.  Note that $\Omega^G_0
=A(G)$, the Burnside ring of $G$, and one can localize at prime
ideals of $A(G)$ as explained in the last chapter of \cite{tD}.
Results along these
lines, at least philosophically related to what we shall do in 
Section \ref{sec:loc} using Kasparov's $KK^G$-theory, may be found
in \cite{Katz1} and \cite{Katz2}.

\section{Localization and $KK$}\label{sec:loc}

We shall rely on the Localization Theorem of Segal (Proposition 4.1 of 
\cite{S2}) as well as its dual formulation for $K$-homology (see for example
\cite{RW1}, Theorem 2.4). We briefly review how this works.

To compute $[D]$ in the $R(G)$-module $K^G_0(M)$, it suffices to compute its
image $[D]_\fp$ in the localizations $K^G_0(M)_\fp$ of $K^G_0(M)$ with respect
to prime (or even just maximal) ideals $\fp$ of $R(G)$. Every such ideal
has a {\em support}, a conjugacy class $(H)$ of cyclic subgroups $H$ of $G$,
with the property that $\fp$ is the inverse image of a prime ideal of $R(H)$
under the restriction map $R(G)\to R(H)$, and such that if
$\fp$ is also the inverse image of a prime ideal of $R(J)$ for some other
subgroup $J$ of $G$, then $H$ is conjugate to a subgroup of $J$.
Let $M^{(H)}$ denote the union of the fixed sets $M^{gHg^{-1}}$
as $g$ runs over the elements of $G$.
The Localization Theorem states:
\begin{thm}[Localization Theorem \cite{S2}]\label{thm:localization1}
Let $M$ be a compact $G$-space with the $G$-homotopy type of a finite
$G$-CW complex.
Let $\fp$ be a prime ideal of $R(G)$ with support $(H)$.
Then the inclusion $M^{(H)}\hookrightarrow M$ induces an isomorphism on
$K$-homology and $K$-cohomology localized at $\fp$.
\end{thm}

We will also need some properties of Kasparov's bivariant $K$-theory
in the equivariant setting, in other words $KK^G$. If $X$ and $Y$ are
{\em locally\/} compact $G$-spaces, then $C_0(X)$ and $C_0(Y)$
(the continuous $\bC$-valued functions vanishing at infinity) are
abelian $C^*$-algebras with $G$-actions, so $KK^G(C_0(X),\,C_0(Y))$
is defined in \cite{K}; for simplicity, we denote this by $KK^G(X,\,Y)$.
This can be extended to a $\bZt$-graded bivariant theory $KK^G_i(X,\,Y)$.
The bivariant groups subsume both $K$-homology and $K$-cohomology
since $KK^G_i(\textrm{point},\,Y)=K_G^{-i}(Y)$, equivariant $K$-theory with
compact supports as defined in \cite{S2}, and $KK^G_i(X,\,\textrm{point})
=K^G_i(X)$, {\em locally finite\/} equivariant $K$-homology for
locally compact spaces. There is an associative bilinear Kasparov product:
\[ \otimes_X\co KK^G_i(Z,\,X)\times KK^G_j(X,\,Y) \to KK^G_{i+j}(Z,\,Y).\]
If $E$ is a $G$-vector bundle over a compact 
$G$-space $X$, then $E$ corresponds to a finitely generated projective
module over $C(X)$ with compatible $G$-action. Since $C(X)$ is
commutative, we may view this as a (Kasparov) $C(X)$-bimodule, which 
gives a class 
\begin{equation}\label{eq:doublebracket}
[[E]]\in KK^G(X,\,X),
\end{equation}
which is in fact the Kasparov product $[E]\otimes_X [\Delta_X]$
of $[E]\in K_G^0(X)$ with the $KK^G$-class of the diagonal map $\Delta_X\co
X\to X\times X$ (see \cite{Bl}, Lemma 24.5.3). The cup-product
in $K^*_G$ may then be expressed in this language since
\[ [E]\cup [E'] = [E]\otimes_X [[E']] = [E']\otimes_X [[E]]. \]
We can also generalize Theorem \ref{thm:localization1} as:
\begin{thm}[Localization in $KK^G$]\label{thm:localization2}
Let $\fp$ be a prime ideal of $R(G)$ with support $(H)$, and let $X$ and
$Y$ be locally compact $G$-spaces which each have the $G$-homotopy type
of $W_1\setminus W_2$, $(W_1,\,W_2)$ some finite $G$-CW pair.
Then the inclusions $X^{(H)}\hookrightarrow X$ and 
$Y^{(H)}\hookrightarrow Y$ induce isomorphisms
\[ KK^G_i(X^{(H)},\,Y^{(H)})_\fp \cong KK^G_i(X,\,Y)_\fp. \]
\end{thm}
\begin{proof}
By the long exact sequences for the pairs $(X,\,X^{(H)})$ and $(Y,\,Y^{(H)})$,
it is enough to show that $KK^G_i(X,\,Y)_\fp=0$ if $X^{(H)}=\emptyset$
or $Y^{(H)}=\emptyset$. Equivariant homotopy invariance, the long exact
sequences, and inductions on cells reduce everything to the case
of a single equivariant cell in each variable, and then by Bott periodicity,
we just need to show that if $K$ and $J$ are subgroups of $G$,
\[ KK^G_*(G/K,\,G/J)_\fp=0 \]
whenever $H$ is not subconjugate to both $K$ and $J$. But by equivariant
Poincar\'e duality (\cite{KInv}, \S4), we may move the 
($0$-dimensional) $G$-manifold $G/K$ across to the other side, obtaining
that
\[ KK^G_*(G/K,\,G/J) \cong K_G^*((G/K) \times (G/J)) .\]
But $(G/K) \times (G/J)$ only has $H$-fixed points if both $G/K$ and $G/J$ do,
and we conclude using Theorem  \ref{thm:localization1} (or the fact
from \cite{S1} on which it is based, that $R(K)_\fp\ne 0$ if and
only if $H$ is subconjugate to $K$).
\end{proof}

Now let's return to the situation of Section \ref{sec:intro}.
For any $G$-invariant open subset
$U$ of $M$, we have a restriction map $K^G_*(M)\to K^G_*(U)$ sending
the class $[D_M]$ to the class $[D_U]$ of the signature operator on $U$
(with respect to some complete $G$-invariant metric on $U$ --- see the
introductions to \cite{Ro2} and \cite{Ro3} for more details and references).
If $U$ is an nice open neighborhood of $M^{(H)}$, the Localization Theorem
again says that the restriction map $K^G_*(M)_\fp\to K^G_*(U)_\fp$ is
an isomorphism, so that $[D_M]_\fp$ may be identified with $[D_U]_\fp$.
Passing to the limit over smaller and smaller $G$-invariant neighborhoods
$U$ of $M^{(H)}$, we obtain:
\begin{thm}\label{thm:localization3}
\textup{(\cite{RW1}, Theorem 2.6. However, the result was not
stated correctly there when $H$ is not normal in $G$; see also
\cite{Ro3}, Theorem 2.9.)}
Let $M$ be an oriented closed $G$-manifold, where $G$ is a finite
group acting smoothly on $M$ and preserving the orientation, and let $H$
be a cyclic subgroup of $G$. Then $[D]_\fp$ is a sum of terms
coming from the various components $F_i$ of $M^{(H)}/G=
M^H/N$, where $N=N_G(H)$ is the normalizer of $H$ in $G$. 
The contribution from $F_i$ only depends on the germ 
of $G\cdot \overline F_i$ in $M$ as a $G$-space.  
\textup{(}Here if $F_i$ is a component 
of $M^H/N$, $\overline F_i$ is its preimage in $M^H$, which might
be disconnected.\textup{)}
\end{thm}

If $H$ is not normal in $G$, then $M^{(H)}$ can fail to be a manifold,
and Theorem \ref{thm:localization3} is of only limited usefulness.
Hereafter we will ignore this situation, and
assume $H$ is normal in $G$. (Even when this is not the case, we can obtain
some useful information by replacing $G$ by the normalizer of $H$
or even something smaller; see Theorem \ref{Gtocyclic} below.)
In  fact, if $G$ is abelian or a quaternion group, then every cyclic subgroup
of $G$ is normal, so we can replace $M^{(H)}$ by the manifold $M^H$ in 
Theorem \ref{thm:localization3}. 
Note that even when $H$ is normal in $G$, $M^H$ may still be
disconnected, and the $G$-action on it may permute the components. However,
if $F$ is a component of $M^H$, and if $G'$ is the (setwise) stabilizer of
$F$ in $G$, then $G\cdot F$ is the disjoint union of $|G/G'|$ components,
and the contribution of $G\cdot F$ to $[D_M]_\fp\in K^G_*(M)_\fp$ 
may be identified
with the class in $K^{G'}_*(U)_\fq$ of the signature operator on some
small $G'$-invariant tubular neighborhood $U$ of $F$, where
$\fq$ is the prime idea of $R(G')$ corresponding to $\fp\triangleleft
R(G)$. (See the beginning of the proof of Theorem 2.12 in \cite{Ro3}.) 
To avoid cluttering up the notation, we thus replace $G$ by $G'$ and
assume that $H$ is a cyclic normal subgroup of $G$, that $F$ is a component
of $M^H$, that $G$ acts on $F$, and that $U$ is a $G$-invariant tubular
neighborhood of $F$. 
We want to compute the class $[D_U]_\fp\in K^G_*(U)_\fp
\cong K^G_*(F)_\fp$ in terms of the class $[D_F]_\fp\in K^G_*(F)_\fp$
(or a twisted analogue, if $F$ is not orientable) and the
normal bundle of $F$. This will require looking at 
the signature operator along the fibers of a vector bundle, in the
case of our specific situation. First it will be convenient to point
out certain facts about ``change of group'' in equivariant $K$-(co)homology.

\begin{thm}[see \cite{Ru}]\label{Gtocyclic}
Let $G$ be a finite group, let $X$ be a locally compact $G$-space, and
let $r\co K^G_*(X)\to \bigoplus_{S\subseteq G\ \mathrm{cyclic}}K^S_*(X)$ be
the direct sum of the restriction maps from $G$-equivariant $K$-homology to
$S$-equivariant $K$-homology, as $S$ runs over the cyclic subgroups of $G$.
Then the kernel of $r$ is torsion of exponent dividing the order of $G$.
If $\fp$ is a prime ideal of $R(G)$ with support $(H)$, then modulo
torsion of exponent dividing the order of $G$, 
$r_\fp\co K^G_*(X)_\fp\to \bigoplus_{H\subseteq
S\ \mathrm{cyclic}}K^S_*(X)_\fp$ is injective.
\end{thm}
\begin{proof}This is proved in \cite{Ru} (in the dual situation of 
$K$-cohomology, for the much harder case of compact Lie groups, but without
the statement about the exponent of the torsion). 
In our particular situation the proof is easy once one makes use of
Artin's Theorem on induced characters (\cite{Serre}, \S II.9.4),
which asserts that for any $\chi\in R(G)$, $|G|\chi$ is an
integral  linear combination
of characters induced from cyclic subgroups. We only need this for $\chi=1$,
the trivial representation. Write $|G|=\sum_S \mathrm{Ind}_S^G\chi_S$,
where $S$ runs over the cyclic subgroups of $G$ and
with $\chi_S\in R(S)$. Then construct a map
\[ s\co \bigoplus_{S\subseteq 
G\ \mathrm{cyclic}}K^S_*(X) \to K^G_*(X)
\]
by sending $c\in K^S_*(X)$ to $\mathrm{Ind}_S^G\left( \chi_S\cdot c\right)$,
where $\mathrm{Ind}_S^G$ denotes the composite
\[ K^S_*(X)\cong K^G_*(G/S\times X) \stackrel{(\mathrm{proj}_2)_*}
{\longlongarrow} K^G_*(X). \]
Then by construction, $s\circ r$ is multiplication by $|G|$, and so the
kernel of $r$ is torsion of exponent dividing $|G|$.

The final statement about the localized case follows now from \cite{S1},
Proposition 3.7, which asserts that $R(S)_\fp=0$ unless a conjugate
of $S$ contains $H$, together with the fact that in the above construction,
we really only needed one cyclic subgroup in each conjugacy class of
cyclic subgroups (since all conjugate subgroups induce the same representations
of $G$).
\end{proof}

An immediate application of Theorem \ref{Gtocyclic} is that, at the expense
of killing some torsion of exponent dividing $|G|$, we can always
restrict attention to cyclic groups, thereby bypassing the problem we
mentioned earlier about $M^{(H)}$ not always being a manifold.

\begin{lem}\label{lem:GtoH}
Let $G$ be a finite abelian group, let $\fp$ be a prime ideal
of $R(G)$, and let $H$ be the cyclic subgroup which is its support.
Then there is a unique prime ideal of $R(H)$, say $\fq$,
which pulls back to $\fp$
under restriction $r\co R(G)\to R(H)$, and $R(H)_\fp = R(H)_\fq$.
Furthermore, if the residual characteristic $p$ of $\fp$ is either $0$ or
relatively prime to $|G/H|$, then $r$ induces an isomorphism $R(G)_\fp
\stackrel{\cong}{\to} R(H)_\fq$. 
\end{lem}
\begin{proof}
Since $G$ is abelian, $R(G)=\bZ\widehat G$ and $R(H)=\bZ\widehat H$,
where $\widehat G$ and $\widehat H$ are the dual groups. Since
$\widehat G\twoheadrightarrow \widehat 
H$, $R(H)\cong R(G)/I$, where $I$ is the
kernel of $r$. By \cite{S1}, Proposition 3.3(i), $N_G(H)/Z_G(H)$ acts
transitively on the prime ideals of $R(H)$ pulling back to $\fp$, so
in the abelian case there is only one such ideal, say $\fq$,
and $\fp/I=\fq$. We have $R(H)_\fp = R(H)_\fq$ by \cite{Bour},
Ch.\ II, \S2.2, Proposition 6. Furthermore,
the character $\chi\in R(G)$ of $\mathrm{Ind}_H^G 1_H$ takes the value 
$|G/H|$ on $H$
and vanishes off of $H$, so that if $p$ is $0$ or relatively prime to 
$|G/H|$, then $\chi\notin \fp$ and $\chi$ annihilates $I$,
so that $I_\fp=0$ (\cite{Bour}, Ch.\ II, \S2.2, Corollary 2)
and $r$ induces an isomorphism $R(G)_\fp \stackrel{\cong}{\to} R(H)_\fq$ 
(\cite{Bour}, Ch.\ II, \S2.5, Proposition 11). 
\end{proof}
\begin{lem}\label{lem:infl}
Let $F$ be a closed oriented $G$-manifold,
with $G$ preserving the orientation,
and suppose that a normal subgroup $N$ of $G$ acts trivially on $F$ 
\textup{(}so that
the $G$-action on $F$ comes from an action of $G/N$\textup{)}. Then the class
of the $G$-equivariant signature operator of $F$ in $K_*^G(F)$ is the 
image of the $(G/N)$-equivariant
class of the signature operator in $K_*^{G/N}(F)$, under the ``inflation''
map of $R(G/N)$-modules $K_*^{G/N}(F)\to K_*^G(F)$.
\end{lem}
\begin{proof}Obvious.\end{proof}
Now we're ready for the key step in the calculation.  First, a few
relevant reminders concerning Kasparov Theory. Suppose $E
\stackrel{p}{\to} B$ is a fibration with smooth manifold fibers,
and we are given a differential operator $D_E$ on $E$ which only involves
differentiation along the fibers, and which is elliptic when
restricted to each fiber. (Note that it is not essential that $E$
itself be a manifold, as the natural domain of $D_E$ consists of
continuous functions on $E$ which are smooth in the fiber directions.)
Then the ``elliptic operator along the fibers'' $D_E$ commutes with
multiplication by functions pulled back from the base $B$ and defines a
Kasparov class in $KK_*(E, B)$, or in the corresponding equivariant
$KK$-group if everything commutes with the action of a finite group.
In fact, this is a special case of what is done in \cite{CoSk},
and is the set-up for proving the index theorem for families
using $KK$-theory. (See also \cite{Hig}, \S4.8.)

Secondly, we need to review something about the calculation of
Kasparov products, as developed for example in \cite{CoSk}, Appendix A,
in \cite{Hig}, \S5, and in \cite{Bl}, \S18. Suppose one has classes
in $KK(E, B)$ and in $KK(B,\bC)$, represented by Kasparov bimodules
$(\cH_1,T_1)$ and $(\cH_2,T_2)$. Thus $\cH_1$ is a $\bZ/2$-graded 
Hilbert $B$-module with an action of $E$, $T$ is an odd $B$-linear operator
on $\cH_1$ ``approximately commuting'' with the action of $E$, and
similarly for the Hilbert space $\cH_2$ and the operator $T_2$.
Then the Kasparov product $[\cH_1,T_1]\otimes_B [\cH_2,T_2]$ is
represented by the $\bZ/2$-graded Hilbert space $\cH=\cH_1\otimes_B 
\cH_2$, together with an operator $T$ which is a $T_2$-connection
in the sense of \cite{CoSk}, Appendix A. In the situation where
$E\stackrel{p}{\to} B$ is a smooth manifold fiber bundle,
$T_1$ comes from an elliptic differential operator $D_E$ of order $1$
along the fibers as above, and $T_2$ comes from
an elliptic differential operator $D_B$ of order $1$ on $B$,
then $T$ comes from the elliptic operator $D_E\otimes 1+1\otimes
p^*D_B$ on $E$.

\begin{thm}\label{thm:productdecomp}
Let $G$ be a finite abelian group, let $\fp$ be a prime ideal
of $R(G)$, and let $H$ be the cyclic subgroup which is its support. 
Let $M^n$ be a smooth
compact $G$-manifold equipped with a $G$-invariant orientation, and
let $[D_M]\in K^G_*(M)$ be the class of the signature operator on $M$.
Observe that $G$ permutes the components $F$ of $M^H$.
Then $[D_M]_\fp\in K^G_*(M)_\fp\cong K^G_*(M^H)_\fp$ is a sum over 
orbits $G\cdot F$ for this action. Let $G'$ be the 
\textup{(}setwise\textup{)} stabilizer in $G$ of some component $F$,
and let $\fq$ be the prime idea of $R(G')$ corresponding to $\fp\triangleleft
R(G)$.  The contribution of $G\cdot F$ to $[D_M]_\fp$ is the
image in $K^G_*(G\cdot F)_\fp\cong K^{G'}_*(F)_\fq$
of the Kasparov product of $[D_F]\in K^{G'/H}_*(F)$, ``inflated''
from $K^{G'/H}_*$ to $K^{G'}_*$ and then localized at $\fq$, 
and of the class in $KK^{G'}(F,F)_\fq$ of the ``signature operator 
along the fibers'' on the normal bundle to $F$ in $M$.
\textup{(}In case $F$ is non-orientable, both the signature operator 
on $F$ and the signature operator along the fibers of the normal
bundle must be twisted by the real line bundle determined by $w_1(F)$.
\textup{)}
\end{thm}
\begin{proof}By the Localization Theorem \ref{thm:localization1},
the contribution of $G\cdot F$ to $[D_M]_\fp$ is the same as
the class of $[D_U]_\fp$, where $U$ is a $G$-invariant tubular
neighborhood of $G\cdot F$, or as
the class of $[D_{U'}]_\fp$, where $U'$ is a $G'$-invariant tubular
neighborhood of $F$. We may identify $U'$ with the total space $E$
of the normal bundle $E \stackrel{p}{\to} F$
to $F$, equipped with a $G'$-invariant metric.
Fixing a $G'$-invariant connection on $E$ enables us to identify
$\twedge \cotan E$ with $p^*(\twedge E^*)\otimes p^*(\twedge \cotan F)$, 
and the signature operator
on $E$ with $D_{\text{fib}}\otimes 1 + 1\otimes p^*D_F$,
where $D_{\text{fib}}$ is the signature operator along the fibers. 
When $F$ is not orientable, then $E$ will not be either, and
both $D_F$ and $D_{\text{fib}}$ have to be taken with coefficients
in the flat real line bundle determined by $w_1(F)$. Also,
it's understood that we have to choose compatible orientations on
the two factors, as discussed in \cite{AS3}, p.\ 581.
So by the above description of the Kasparov product, $[D_E]$ is
the Kasparov product of $[D_{\text{fib}}]$ and $[D_F]$.
(Note: Since $G$ preserves the orientation of
$M$, $F$ has to have even codimension. Thus $[D_{\text{fib}}]$
lives in $KK^{G'}_0$, not $KK^{G'}_1$. This is important since the 
signature operator class on a product is the Kasparov product of the 
signature classes on the factors
\emph{provided} that the manifolds aren't both odd-dimensional
\cite{RW2}; in the exceptional case there is a factor of 2 because of
the way one keeps track of the gradings on the Clifford algebras.
Since the fibers of $E$ have even dimension we don't have a problem here.)
Finally, since $H$ acts trivially on $F$, 
we may apply Lemma \ref{lem:infl} to view $F$ as a ($G'/H$)-manifold.
\end{proof}
As long as $H$ is non-trivial and acts effectively, 
Theorem \ref{thm:productdecomp} reduces the
calculation of $[D_M]_\fp$ down to the situation of smaller manifolds
$F$ and smaller groups $G'/H$, provided we can compute the
contribution of the signature operator along the fibers. We do this
calculation next. Here the manifold structure of $F$ becomes
irrelevant, and all we care about is the $G$-vector bundle $E$.
(Note that for simplicity of notation, we have converted $G'$ to $G$.)
\begin{prop}\label{prop:fiberoperator}
Let $F$ be a compact $G$-space, where $G$ is a finite group, and let
$\fp$ be a prime ideal of $R(G)$ with support a central cyclic subgroup $H
=\langle g_0\rangle$
that acts trivially on $F$. Let $E$ be a real
$G$-vector bundle over $F$ of even dimension $2k$, and assume that
$E$ can be given a $G$-invariant orientation. Note that
$E$ splits as a direct sum of isotypical subbundles for the various
irreducible real representations of $H$. Assume that the trivial
representation of $H$ doesn't occur in this decomposition. Since $H$ is cyclic,
the remaining irreducible representations of $H$  are all two-dimensional 
and of complex type, with the exception of the ``sign'' representation
$G\twoheadrightarrow \{ \pm 1\}$ if $|H|$ is even. So
we may write $E$ as a direct sum of oriented even dimensional subbundles
$E(-1)$ and $E(e^{i\theta_j})$, where $0<\theta_j< \pi$. Here $g_0$ acts
on $E(-1)$ by multiplication by $-1$\textup{;} for $0<\theta_j< \pi$,
$E(e^{i\theta_j})$ has a complex structure and $g_0$ acts on it by complex
multiplication by $e^{i\theta_j}$. Let $E_c(e^{i\theta_j})$ denote
$E(e^{i\theta_j})$ viewed as a complex vector bundle. Since $H$ is central, 
the decomposition of $E$ into the $E(e^{i\theta_j})$'s is preserved by $G$.
Furthermore, since the complex structure on  $E_c(e^{i\theta_j})$ comes
from the action of $g_0$, this is a complex $G$-vector bundle and not
just a complex $H$-vector bundle.  Let $D^{\mathrm{fiber}}_E$ be the
signature operator along the fibers of $E$, for a $G$-invariant Euclidean
metric on $E$. \textup{(}So on each fiber of $E$, $D$ looks like the signature
operator on $\bC^k$, where $g_0$ acts on $\bC^k$ with eigenvalues 
$-1$ and/or $e^{i\theta_j}$, $0<\theta_j< \pi$.\textup{)}
Then $D^{\mathrm{fiber}}_E$ defines a class
$[D^{\mathrm{fiber}}_E] \in KK^G(E,\,F)$. Localized at $\fp$, this class
lives in $KK^G(E,\,F)_\fp \cong KK^G(F,\,F)_\fp$,
and may be identified with $[[\mathcal{E}]]=[\mathcal{E}]\otimes_F [\Delta_F]$
\textup{(}notation of equation \textup{(\ref{eq:doublebracket}))},
where $\mathcal{E}$ is the cup product of classes 
$[\mathcal{E}(e^{i\theta_j})]$,
$0<\theta_j< \pi$, and $[\mathcal{E}(-1)]$. Furthermore, we have
\[
[\mathcal{E}(e^{i\theta_j})]=\left[\twedge E_c(e^{i\theta_j}) 
\right] / \left(\left[\twedge^{\mathrm{even}}
	E_c(e^{i\theta_j})\right]
	- \left[\twedge^{\mathrm{odd}}
	E_c(e^{i\theta_j})\right]\right) 
\]
for $\theta_j< \pi$, and also for $\theta_j=\pi$ if $E(-1)$ has a
$G$-invariant complex structure. If $E(-1)$ has a
$G$-invariant spin$^c$ structure, we have a similar formula:
\[
[\mathcal{E}(-1)] =\left[ \mathcal{S}\left(E(-1)\right)\right] /
\left(\left[ \mathcal{S}^+\!\left(E(-1)\right)\right]-
      \left[ \mathcal{S}^-\!\left(E(-1)\right)\right]\right),
\]
where $\mathcal{S}\left(E(-1)\right)$ is the complex spinor bundle for
the  spin$^c$ structure on $E(-1)$, and $\mathcal{S}^\pm$ are the
half-spinor bundles. Finally, if $E(-1)$ does not have a $G$-invariant 
spin$^c$ structure, the formula is the same, but must be interpreted
in the sense of twisted coefficients.
\end{prop}
\begin{proof}
Let $\mathcal{X}$ denote the continuous field of Hilbert spaces over
$F$, whose fiber over $x\in F$ is $L^2\left( \twedge \cotan _\bC E_x \right)$,
with the $\bZt$-grading of \cite{AS3}, p.\ 575.
Since $D^{\mathrm{fiber}}_E$ commutes with
multiplication by functions on the base $F$, and is self-adjoint, local, and
elliptic along the fibers, and since $G$ acts isometrically,
the pair $(\mathcal{X},\, D^{\mathrm{fiber}}_E)$
satisfies the conditions for an unbounded Kasparov module in the sense
of \cite{BJ}, giving a class $[D^{\mathrm{fiber}}_E] \in KK^G(E,\,F)$.
Since $+1$ is not an eigenvalue of the action of $g_0$ on $E$, $E^{(H)}=F$
and we have an isomorphism 
\[
KK^G(E,\,F)_\fp \cong KK^G(F,\,F)_\fp
\]
by Theorem \ref{thm:localization2}. But on an even-dimensional 
$\mathrm{spin}^c$ manifold, the signature operator may be expressed as
the Dirac operator with coefficients in the dual of the complex spinor
bundle (\cite{LM}, \S II.6.2). Suppose first that $E(-1)$ also has a 
$G$-invariant complex structure, so that $E$ is the underlying real
$G$-vector bundle of a complex $G$-vector bundle $E_c$,
the direct sum of $E_c(-1)$ and the $E_c(e^{i\theta_j})$'s. Then 
\[
\Cl E^* \cong \twedge(E_c^*)\otimes_\bC \twedge(E_c)
\cong \left( \bigotimes_{0 < \theta_j \le \pi} 
\twedge E_c(e^{i\theta_j}) \right) \otimes 
\left( \bigotimes_{0 < \theta_j \le \pi} \twedge
E_c(e^{i\theta_j})^* \right),
\]
with the first factor
identified as the complex spinor bundle and the second factor identified
with its dual. So $[D^{\mathrm{fiber}}_E]$ is the class of the Dirac
operator along the fibers, with coefficients in $\twedge(E_c)$. But the
Dirac operator along the fibers gives the inverse of the Thom isomorphism
$\tau\in KK^G(F,\,E)$ in equivariant $K$-theory, so that we have the 
formula
\[ \tau\otimes_E [D^{\mathrm{fiber}}_E] = [[\twedge(E_c)]]\in KK^G(F,\,F).
\]
When we localize at $\fp$ and restrict to the fixed-point set (which is just
the zero-section of $E$), $\tau$ is multiplication by $\twedge_{-1}(E_c)
=[\twedge^{\text{even}}(E_c)]-[\twedge^{\mathrm{odd}}(E_c)]$
(\cite{S2}, \S3). So
\[ [D^{\text{fiber}}_E] = [[\twedge(E_c) / \twedge_{-1}(E_c) ]] ,
\]
which in turn splits into pieces corresponding to the various
rotation angles $\theta_j$, 
as claimed.  (The division makes sense since $\twedge_{-1}(E_c)$ is a unit
in equivariant $K$-theory localized at $\fp$.)

Now consider the case where $E(-1)$ does not have a complex structure (or
at least one that is $G$-invariant). We proceed as before, except that
we need twisted coefficients in the Thom isomorphism (\cite{K}, \S5, Theorem
8) for the factor associated to $E(-1)$. Note that in the case of a
$G$-invariant spin$^c$ structure on $E(-1)$, $\mathcal{S}^+\!
\left(E(-1)\right)$ substitutes for $\twedge^{\text{even}}E_c(-1)$,
and $\mathcal{S}^-\!
\left(E(-1)\right)$ substitutes for $\twedge^{\text{odd}}E_c(-1)$.
\end{proof}
\begin{exm}To give a very simple example, suppose $G$ is abelian,
$\fp$ is a prime ideal with support $H$, and $F$ is a component
of $M^H$ whose normal bundle is stably equivariantly trivial
(i.e., a $G$-invariant tubular neighborhood of $F$
is stably just a product of $F$ with a representation
space $V$ of $G$ whose restriction to $H$ is non-trivial on a
generator  of $H$). As we've seen, $V$ must have a $G$-invariant
complex structure, so we can think of $V$ as the realification of
a complex representation $V_c$, and we choose the orientation of $F$
so that the orientation on $F\times V_c$ agrees with the orientation
of $M$. Then our formula for the contribution of $G$ to $[D_M]_\fp$
reduces simply to
\[ \left[\twedge (V_c)/\twedge_{-1} (V_c)\right]\cdot [D_F].\]
Here $\twedge (V_c)$ and $\twedge_{-1} (V_c)$ are viewed as elements of 
$R(G)$. While we can't divide them in $R(G)$, the fact that each
constituent of $V_c$ is non-trivial on a generator of $H$ means that
$\twedge_{-1} (V_c)$ does not lie in the prime ideal $\fp$, so the
division makes sense in $R(G)_\fp$, the coefficient ring for the
localized theory. To see this, first observe that if we write
$V_c$ as a sum of irreducible characters $\chi_i\co G\to U(1)$, then
$\twedge_{-1} (V_c)=\prod_i(1-\chi_i)$. So we just need to show
that if $\chi$ is a one-dimensional representation of $G$ which is
non-trivial on a generator $g$ of $H$, then $1-\chi\notin \fp$.
If the residual characteristic of $\fp$ is $0$, then $\fp$ is just
the ideal consisting of virtual representations of $G$ whose
characters vanish at $g$, so by assumption on $\chi$, $\chi(g)\ne 1$, i.e.,
$1-\chi\notin \fp$. And if the residual characteristic of $\fp$
is finite, say $p$, then the order of $H$ is relatively prime to
$p$ (\cite{S1}, Proposition 3.5), which since $\chi(g)\ne 1$
will force $\chi|_H$ to map to an element
other than 1 in the finite field $R(H)/\fq$ of characteristic $p$.
(Here $\fq$ is as in Lemma \ref{lem:GtoH}.) Thus again
$1-\chi\notin \fp$ in this case.\qed
\end{exm}
We put everything together to show that the equivariant $K$-homology
class of the signature operator can contain quite complicated
torsion information not preserved under $G$-homotopy equivalences, 
even when $M$ is a very simple manifold and $G$
is cyclic of prime order.
\begin{thm}Let $p$ be an odd prime. Then there exist actions of $G=\bZ/p$
on odd spheres $M=S^{2k+1}$ for which the equivariant $K$-homology
class $[D_M]\in K^G_1(M)$ contains ``arbitrarily complicated
torsion information.''  More precisely, let $q$ be a prime
\textup{(}which may or may not be equal to $p$\textup{)}. Then there
exist actions of $G=\bZ/p$ on spheres $M_1=M_2=S^{2k+1}$ and
$G$-homotopy equivalences $M_1\stackrel{h}{\to} M_2$ such that
$[D_{M_2}]-h_*([D_{M_1}])$ is torsion of order as large a 
power of $q$ as one wants.
\end{thm}
\begin{proof}(\emph{Sketch})
First suppose $q=p$. Then take $M_1$ and $M_2$ to be
free linear $G$-spheres such that the lens spaces $L_1=M_1/G$ and
$L_2=M_2/G$ are homotopy equivalent but not diffeomorphic. Note
that under the $\bZ$-module isomorphism $K_*^G(M_j)\cong K_*(L_j)$,
$[D_{M_j}]$ corresponds to $[D_{L_j}]\in K_1(L_j)$. This class
is computable (see \cite{RW2}, \S2) and is not homotopy-invariant.
(The idea of the calculation is to write the signature operator as 
a Dirac operator with coefficients in the dual of the spinor bundle 
of the cotangent bundle, as in the proof of Proposition
\ref{prop:fiberoperator}.)  In fact, given $r\ge 1$, we can choose
$k$ sufficiently large (depending on $p$ and $r$) so that there
is a homotopy equivalence $h\co L_1\to L_2$ between lens spaces
of dimension $2k+1$ for which $[D_{L_2}]-h_*([D_{L_1}])$ has
order $p^r$ in $K_1(L_2)$, giving us the example we want when we
pull back to the universal covers.

Next, suppose $q\ne p$, and this time consider lens spaces 
$L_1$ and $L_2$ of dimension ${2j+1}$, each with fundamental group 
$\bZ/q^r$, which are homotopy equivalent
but have different signature operator classes, just as above.  
Then $L_1$ and $L_2$ are $\bZ/p$-homology spheres with rationally 
trivial stable normal bundles, so by ``converse Smith theory'' it is known
that they can be realized fixed sets of semifree actions of $G$ on 
spheres $M_1=M_2=S^{2k+1}$, provided that $j<k$ are in an appropriate 
range ($k$ roughly equal to $jp$). (See \cite{Jones}, \S5, and 
\cite{Wgrpact}, \S6.) In this way one can get an equivariant
homotopy equivalence $h\co M_1\to M_2$. We can then compute
$[D_{M_2}]-h_*([D_{M_1}])$, localized at a prime ideal $\fq$ of residual
characteristic $q$ supported on all of $G$, 
using Theorem \ref{thm:productdecomp}. We
obtain a difference between $[D_{L_2}]$, twisted by some normal
characteristic classes, and something similar for $L_1$, transported
over to $L_2$ via $h$.  One can arrange for $[D_{M_2}]-h_*([D_{M_1}])$
to involve large $q$-primary torsion.
\end{proof}

\section{Comparison with the Atiyah-Singer Theorem}\label{sec:AS}
Let us check that the formula for the localization of
$[D_M]$ derived in Section \ref{sec:loc} agrees with the Atiyah-Singer
$G$-Signature Theorem (\cite{AS3}, Theorem 6.12) in the case when
the dimension of $M$ is even. To see this, let $g\in G$ and let
$H$ be the cyclic subgroup it generates. Let $\fp$ be the prime
ideal of $R(G)$ consisting of virtual representations whose characters
vanish at $g$, so that $R(G)/\fp\hookrightarrow \bC$ via evaluation
of characters at $g$. Clearly $\fp$ has support $H$ in the sense of
\cite{S1}. We may compute Sign($g$, $X$) (in the sense of \cite{AS3})
by mapping $[D_M]_\fp$ to $R(G)_\fp$ via the collapse map $c\co M\to 
\mathrm{point}$, and then mapping to the residue field
$R(G)_\fp/\fp_\fp$ (a subfield of $\bC$, in fact a number field). 
We get a sum of terms coming
from the components $F$ of $M^H$, and for purposes of computing
Sign($g$, $X$), we may as well assume $G=H$. Then Theorem 
\ref{thm:productdecomp} and Proposition \ref{prop:fiberoperator}
apply to this situation. So the contribution of $F$ to $[D_M]_\fp$
is thus effectively the Kasparov product of the non-equivariant
$K$-homology class $[D_F]$ of the signature operator on $F$
and of the class in $KK^H(F,F)$ of the signature operator on the fibers
of the normal bundle. The Chern character of the former is the
Poincar\'e dual of the Atiyah-Singer $L$-class $\mathcal{L}(F)$,
or in other words the factor $A_1$ in \cite{AS3}, p.\ 581,
and of the latter is the Chern character of the virtual $H$-bundle
described in Proposition \ref{prop:fiberoperator}. The second
factor accounts for the factors $B_1$ and $C_1^\theta$ in \cite{AS3},
p.\ 581. Consider for example the contribution of $E_c(e^{i\theta})$,
which in the Atiyah-Singer notation is $N^g(\theta)$,
when $0 < \theta < \pi$. If this splits into complex line bundles
with first Chern classes $x_j$, then since $g$ acts by $e^{in\theta}$
on the $n$-th tensor power of one of these line bundles, we get 
from $[[\twedge(E_c(e^{i\theta})) / \twedge_{-1}(E_c(e^{i\theta}))]]$ a 
contribution of
\[
\prod_j \left(\frac{1+e^{x_j}e^{i\theta}}
{1-e^{x_j}e^{i\theta}}\right)
=\prod_j 
\left(\frac{e^{(x_j+i\theta)/2}+e^{-(x_j+i\theta)/2}}
{e^{-(x_j+i\theta)/2}-e^{(x_j+i\theta)/2}}\right)
,
\]
which up to a sign ($(-1)^{s(\theta)}$) is
\[
\prod_j \coth \left(\frac{x_j+i\theta}{2}\right),
\]
just as on p.\ 581 of \cite{AS3}.

\section{Extension to the non-smooth case,\\ Concluding
remarks}\label{sec:nonsmooth}
The results of Section \ref{sec:loc} can be interpreted
as an inductive algorithm for computing the class 
in equivariant $K$-homology of the signature operator $D_M$
of a smooth closed $G$-manifold $M$, modulo perhaps the loss of 
some torsion of order dividing a power of the order of $G$, on 
the basis of two ingredients:
\begin{enumerate}
\item the (non-equivariant) $K$-homology classes of the signature
operators on certain submanifolds, namely, the connected components $F$
of the fixed sets for cyclic subgroups $H\subseteq G$. (In case $F$
is non-orientable, we use the signature operator with coefficients
in the real line bundle determined by $w_1(F)$.)
\item certain characteristic classes in equivariant $K$-theory 
for the normal bundles $E$ of these submanifolds $F$, as given in
Proposition \ref{prop:fiberoperator}. (If $F$ is not orientable,
then the normal bundle isn't orientable either, and we replace it
by its tensor product with the real line bundle determined by
$w_1(F)$, which now \emph{is} orientable.)
\end{enumerate}
Before going on to the non-smooth
case, let us review this algorithm. By Theorem \ref{Gtocyclic}, if
we are prepared to accept the loss of some torsion of order dividing
the order of $G$, we can always restrict to subgroups and reduce to 
the case where $G$ is abelian, in fact cyclic. Then since 
\[
  K_*^G(M) \to \bigoplus_{\fp\text{ maximal in }R(G)} K_*^G(M)_\fp
\]
is injective (\cite{Bour}, Ch.\ II, \S3.3, Theorem 1), it is no 
loss of generality to localize at a maximal ideal $\fp$ of $R(G)$, say 
with support $H$. If $H=\{1\}$, then $M^H=M$ and localization
doesn't do much; however, if the residual characteristic $p$ of
$\fp$ is prime to $|G|$, then by Lemma \ref{lem:GtoH}, $R(G)_\fp
\cong \bZ_{(p)}$ and we loose nothing by forgetting $G$ entirely.
If $H=\{1\}$ and the residual characteristic $p$ of
$\fp$ does divide $|G|$, then $R(G)_\fp \to \bZ_{(p)}$ is not
an isomorphism, but the restriction map $K_*^G(M)_\fp\to 
K_*(M)_{(p)}$ only kills some $p$-primary torsion.
So assume $H\ne \{1\}$. If $H$ fixes all of $M$, then we can apply
Lemma \ref{lem:infl} (with $N=H$), and deduce that the class
$[D_M]\in K_*^G(M)_\fp$ comes from the class $[D_M]\in K_*^{G/H}(M)$,
which only involves the action of the smaller group $G/H$. Otherwise,
choose some component $F$ of $M^H$, which is now a submanifold
of smaller dimension, and apply Theorem \ref{thm:productdecomp}.
This computes the contribution of $F$ to $[D_M]_\fp$ in terms of
the ingredients (1) and (2) above.

There is hope for carrying out all or most of the same program when
$M$ is only a Lipschitz manifold and the action of $G$ is Lipschitz 
and locally linear, using the Lipschitz signature operator and its
$KK$-class as constructed in \cite{Tel2}, \cite{Tel3}, \cite{Hil1}, 
and  \cite{Hil2}. The Lipschitz locally linear category of group 
actions was studied
to some extent in \cite{RW1} and in \cite{RotW}, and as explained
in \cite{RotW}, is quite close to the topological locally linear
category. (Since PL manifolds have a canonical Lipschitz
structure, the discussion here includes the PL locally linear case.
However, the construction of the PL signature operator class in
\cite{Tel1} is much easier than in the Lipschitz case.)

The first steps of the program, involving restriction to cyclic
subgroups and localization at prime ideals of $R(G)$, go through
with almost no change, thanks to \cite{RW1}, which enables us to
localize the Lipschitz signature operator in a $G$-invariant
neighborhood of some fixed set component $F$. The problem is
that even in the PL locally linear category, this neighborhood
can be identified with a block bundle over $F$, not in general 
with a vector bundle, and it is not clear if one can split the 
signature operator as in Theorem \ref{thm:productdecomp}. So we
conclude with the following question:
\begin{ques}Suppose $M$ is a locally linear Lipschitz $G$-manifold,
equipped with a $G$-invariant orientation,
with $G$ a finite abelian group, and let $[D_M]$ be the class in
$K_*^G(M)$ of the Lipschitz signature operator on $M$. Let 
$\fp$ be a prime ideal of $R(G)$,
and let $H$ be the cyclic subgroup which is its support. Let $F$
be a component of $M^H$ \textup{(}a locally flatly embedded topological
submanifold of $M$\textup{)}. Then is it possible as in Theorem 
\ref{thm:productdecomp} to split the contribution
of $G\cdot F$ to $[D_M]_\fp$ as the Kasparov product of $[D_F]$,
the signature operator class on $F$, and of a term corresponding to
the ``signature operator along the fibers''?  If so, how can one
compute the latter?
\end{ques}

\bibliographystyle{amsplain}

\end{document}